\documentclass[a4paper,10pt]{amsart}

\usepackage{amsmath}
\usepackage{amssymb}
\usepackage{graphicx}

\usepackage{amsthm}
\theoremstyle{definition}
\newtheorem*{dfn}{Definition}
\newtheorem*{exa1}{Example 1}
\newtheorem*{exa2}{Example 2}
\newtheorem*{exa3}{Example 3}

\newtheorem*{thm1}{Main Theorem}

\newtheorem*{pro1}{Proposition  1}
\newtheorem*{pro2}{Proposition  2}
\newtheorem*{pro3}{Proposition  3}
\newtheorem*{pro4}{Proposition  4}
\newtheorem*{lem}{Lemma}

\newtheorem*{remark1}{Remark  1}

\newtheorem*{remark4}{Remark4 4}

\newtheorem*{remark5}{Final remark on the main theorem}

\newtheorem*{question1}{Question   1}
\newtheorem*{question2}{Question   2}
\newtheorem*{question3}{Question   3}
\newtheorem*{question4}{Question   4}
\newtheorem*{question5}{Question   5       (A weak version of the Kaplansky-Kadison  conjecture)}

\theoremstyle{plain}

\newtheorem*{coro2}{Corallary 2      (Borsuk-Ulam Theorem)  }
\newtheorem*{coro1}{Corallary 1   }
\newtheorem*{coro3}{Corallary 3   }

\begin{document}

\title[BORSUK-ULAM THEOREM]{A BANACH ALGEBRAIC APPROACH TO THE BORSUK-ULAM THEOREM}
\author{Ali Taghavi}

\address{Faculty of Mathematics and Computer Science,  Damghan  University,  Damghan,  Iran.}
\email{taghavi@du.ac.ir}

%\date{April 14, 2011}

\subjclass [2000]{46S60, 46L85}

\keywords{Borsuk-Ulam theorem ,Graded Banach algebras}

\begin{abstract}
Using methods from the theory of commutative graded Banach algebras, we obtain a generalization of the two dimensional Borsuk-Ulam theorem as follows: Let  $\phi:S^{2} \rightarrow S^{2}$ be a homeomorphism of order n and $\lambda\neq 1$ be an nth root of the unity,   then for every complex valued continuous function $f$ on $S^{2}$,  the function $\sum_{i=0}^{n-1} \lambda^{i}f(\phi^{i}(x))$ must  vanish at some point of $S^{2}$. We also discuss about some noncommutative versions of the Borsuk- Ulam theorem.
\end{abstract}

\dedicatory{ To my children Tarannom and Pouya.}

\maketitle

\section*{Introduction}
\noindent
The classical Borsuk-Ulam theorem states that for every continuous function \;\;\;\;\   $g: S^{n}\rightarrow \mathbb{R}^{n}$,  there always
exist a point $x\in S^{n}$ such that $g(-x)=g(x)$. If we define $f(x)=g(x)-g(-x)$, we obtain an equivalent statement as follows: For every odd continuous function f : $S^{n}\rightarrow \mathbb{R}^{n}$, there exist a point $x\in S^{n}$ such that $f(x)=0$.\\
 We consider the case $n=2$ and identify $\mathbb{R}^{2}$ with the complex numbers $\mathbb{C}$. Let $A=C(S^{2})$ be the Banach algebra of all continuous complex valued functions on $S^{2}$ with the $Z_{2}$-graded structure
 \begin{equation*}
 A=A_{ev}\bigoplus A_{odd}
  \end{equation*}
  where $A_{ev}$ is the space of all even functions and $A_{odd}$ is the space of  all odd functions and the decomposition is the standard decomposition of functions  to even and odd functions.  Then the two  dimensional Borsuk-Ulam theorem says that a homogenous  element of $A$  with   non zero degree, namely an odd function,  is not invertible.\\

   In this paper we are mainly interested in invertible  elements of a graded unital  Banach algebra which are homogenous of nontrivial degree. Some natural questions about   such elements are as follows: If any such element is invertible, can  it be connected to the identity in the space of  invertible  elements?
 What can be said about the relative position of their spectrum with respect to the origin?\\

   As we will see in the main theorem of this paper, for a commutative Banach algebra without nontrivial idempotent, which is graded by a finite Abelian group,
a nontrivial homogenous element cannot be connected to the identity.
    On the other hand, using an standard lifting lemma in the theory of covering spaces, we conclude that an invertible element of $C(S^{2})$ has a logarithm then it lies  in the
same connected component of the identity. This shows that an odd element of $C(S^{2})$ cannot be invertible, so it would gives us a Banach algebraic proof of the Borsuk-Ulam theorem, in dimension two. The purpose of this paper is to translate the classical Borsuk-Ulam theorem into the language of noncommutative geometry.\\

 We also give a concrete example of an $S_{3}$-graded structure for $C_{red}^{*}(F_{2})$, the reduced $C^*$ algebra of the free group on two generators, such that  a nontrivial homogenous element lies in the same component of the identity. This would shows that the commutativity of grading group and graded algebra are necessary conditions in our main theorem. Finally we give a question as a weak version of the Kaplansky-Kadison conjecture. This question naturally arise from our main result, see question 5 at the end  of the paper.\\

\section{Preliminaries}

Let $ A$ be a unital complex Banach algebra and $G$ be a finite  group with neutral  element $e$. A G-graded structure for $A$ is a decomposition
\begin{equation*}
A=\bigoplus_{g\in G} A_{g}
\end{equation*}
 where each $A_{g}$ is a Banach subspace of $A$  and $A_{g}A_{h}\subseteq A_{gh}$. An element $a\in A_{g}  , g\in G $ is called a homogenous element, it is called nontrivial homogenous if $a\in A_{g}$ where $g\neq e$. When $G=Z_{2}$ and $A=A_{0}\bigoplus A_{1}$, an element of $A_{1}$ is called an odd element.  A morphism   $\alpha: A \rightarrow A$ is called a graded morphism provided that $\alpha(A_{g})\subseteq A_{g}$ for all $g \in G$.  Let $A$ be a $G$-graded Banach algebra and $H$ be a normal subgroup of $G$. Then there is an induced $G/H$-graded structure for $A$ with
 \begin{equation*}
 A_{g/H}=\bigoplus_{h\in H} A_{gh}
 \end{equation*}
 The following proposition is used in the proof of our main theorem:
 \begin{pro1}
  Let $G$ be Abelian and $a$ be a nontrivial homogenous element of a $G$- graded Banach algebra, then for some positive integer $n$, there is a $Z_{n}$-graded structure for $A$ such that $a$ is a nontrivial homogenous element of $A$ as a $Z_{n}$-graded algebra.
  \end{pro1}

  \begin{proof}
  This can be proved by induction on order of $G$ as follows: The first step of the induction is obvious sine the only group of order 2 is $Z_{2}$. Let $G$ be a finite Abelian group and $g\in G$ where $g \neq e$. If $G$ is not a cyclic group, then there is a subgroup $H$ which does not contain $g$, so $a$ is a nontrivial homogenous element of the induced $G/H$-graded structure for $A$.  Now the order of $G/H$ is strictly less than the order of $G$. So an induction argument on order of $G$ completes the proof.
   \end{proof}
   Note that the existence of a $Z_{n}$-graded structure for a Banach algebra $A$ is equivalent to existence of a bounded
multiplicative operator $T : A \rightarrow A$ with $T^{n}=Id$. For any such operator, we choose a root of unity $\lambda\neq 1$ and observe that the decomposition
\begin{equation*}
A=\bigoplus_{i=0}^{n-1} \ker (T-\lambda^{i})
\end{equation*}
is a $Z_{n}$-graded  structure. Conversely for the grading $A=\bigoplus_{i=0}^{n-1} A_{i}$,  the following multiplicative operator $T$ satisfies $T^{n}=Id$
\begin{equation*}
T(\sum_{i=0}^{n-1} a_{i})=\sum_{i=0}^{n-1}\lambda^{i}a_{i}.
\end{equation*}
\\

\begin{exa1}
Let $X$ be a compact topological space and $\phi$ be an order $n$ homeomorphism of $X$, that is $\phi^{n}=Id$. Define
$T:C(X) \rightarrow C(X)$ with $T(f)=f\circ \phi$,  where $C(X)$ is the space of all continuous functions on $X$. Then $T$ is a bounded multiplicative operator on Banach algebra $A=C(X)$ with $T^{n}=Id$, so we would have a $Z_{n}$-graded structure for $A$. As a particular example, the  rotation of circle by $2\pi/n$ is an order $n$ homeomorphism of circle. Then  we naturally obtain a $Z_{n}$-graded structure for $C(S^{1})$.
\end{exa1}

For a group $G$, a $Z_{n}$-partition for $G$ is a partition of $G$ into disjoint subsets $G_{i}$, $i=0,1, \ldots n-1$ such that  $G_{i}G_{j}\subseteq G_{i+j\; (mod\; n)}$. This is equivalent to say that  $G_{0}$ is a normal subgroup of $G$ whose quotient is isomorphic to $Z_{n}$. This is an obvious consequence of the normal property for subgroups, see  \cite{HERS}. So every group $G$ which has a normal subgroup $H$ such that $G/H$ is isomorphic to $Z_{n}$ posses a $Z_{n}$-graded structure. For example, for every group $K$, the group  $G=K\times Z_{n}$ has $K\times\{0\}$ as a normal subgroup and $G/K\times\{0\}$ is isomorphic to $Z_{n}$. In this case the $Z_{n}$ partition for $G$ is $G=\bigsqcup _{i\in Z_{n}} K\times\{i\}$  \\

For a discrete group $G$, we denote by $\mathbb{C}G $, the group algebra of $G$ with complex coefficients, namely the space of all linear combinations  $\sum a_{g}g$ where $a_{g}$'s are complex numbers. let $l^{2}(G)$ be the Hilbert space of  all $\gamma:G \rightarrow \mathbb{C}$ such that $|\gamma|_{2}^{2}=\sum_{g\in G} |\gamma(g)|^{2}< \infty$. $G$ acts on $l^{2}(G)$ with $g.\gamma (h)=\gamma(g^{-1}h)$.  This  action defines a unitary representation of $G$ on $l^{2}(G)$. We extend this action by linearity to $\mathbb{C}G $. So each element of $\mathbb{C}G$ can be considered as an element of $B(l^{2}(G))$, the space of all bounded operators on $l^{2}(G)$ and every element of $G$ can be considered as a unitary operator on $l^{2}(G)$.
The reduced C* algebra of $G$, $C_{red}^{*}(G)$, is the closure  of $\mathbb{C}G $,  with respect to operator norm
$|\;  . \;|_{op}$ defined on $B(l^{2}(G))$.\\
 In proposition 1, we shall   prove that a $Z_{n}$-partition for a Group $G$, gives a $Z_{n}$-graded structure for $C_{red}^{*}(G)$. For the proof of the proposition 1, we need to the following technical lemma:
\begin{lem}
Let $(E, |\;.\;|)$ be  a Banach space with a dense linear subspace F. Assume $F=\bigoplus_{i=0}^{n} F_{i}$ is a direct sum of its linear subspaces $F_{i}$'s. Suppose that
the original norm of $F$ is equivalent to the the direct limit norm $|\sum f_{i}|=\sum |f_{i}|$. Then $E=\bigoplus_{i=0}^{n} \overline{F_{i}}$
\end{lem}

\textbf{Proof of lemma}\;\;\; Let  $\sum_{i=0}^{n} \widetilde{f^{i}}=0$ where $\widetilde{f^{i}}$ is in $\overline{F_{i}}$, the closure of $F_{i}$. Then there are sequences $\{f_{k} ^{i}\}_{k\in\mathbb{N}}$ of elements of $F_{i}$  which converges  to $\widetilde{f^{i}}$. So  $ \lim_{k\longrightarrow \infty}$    $|\sum_{i=0}^{n} f_{k}^{i}|=0$. Since the original norm of $F$ is equivalent to the direct sum norm  we conclude that each sequence
$\{f_{k} ^{i}\}$ converges to zero. So $\widetilde{f^{i}}=0$ for all i. From continuity  of norm we have the equivalency of norm and direct sum norm on the space $\bigoplus_{i=0}^{n} \overline{F_{i}}$. So this direct sum is a topological direct sum. let $x\in E$ is given, there is a sequence $\sum _{i=0}^{n}f^{k}_{i}$,  $k\in \mathbb{N}$,  which converges to $x$. This shows that each sequence $\{f^{k}_{i}\}$ is a cauchy sequence which converges to an element $\widetilde{f_{i}}\in \overline{F_{i}}$. So $x=\sum_{i=0}^{n} \widetilde{f_{i}}$.   This completes the proof of lemma.

\begin{pro2}
A $Z_{n}$-partition structure for a group $G=\coprod_{i=0}^{n-1} G_{i}$ gives a $Z_{n}$-graded structure for
$C_{red}^{*}(G)=\bigoplus_{i\in Z_{n}} C_{red}^{*}(G_{i})$.
\end{pro2}

 \begin{proof}
 Assume that $G=\coprod_{i\in Z_{n}} G_{i}$ is a $Z_{n}$-partition for $G$, then $\mathbb{C}G=\bigoplus \mathbb{C}G_{i}$ and $l^{2}(G)=\bigoplus l^{2}(G_{i})$  where
\begin{equation*}
l^{2}(G_{i})=\{\gamma\in l^{2}(G):\gamma|_{G_{j}}=0
 \;\; j\neq i\}.
   \end{equation*}
    For $x_{i}\in \mathbb{C}G_{i}$ and $\gamma_{j}\in l^{2}(G_{j})$,   we have $x_{i}.\gamma_{j}$ belongs to $l^{2}(G_{j+i \;\;(mod\; n) })$.  Assume $x=\sum_{i=0}^{n-1}  x_{i}  $ and $\gamma=\sum_{i=0}^{n-1} \gamma_{i}$. Then
   \begin{equation*}
   |x.\gamma_{j}|_{2}^{2} =\sum_{i=0}^{n-1} |x_{i}.\gamma_{j}|_{2}^{2}\;\;\;and \;\;|x_{i}.\gamma|_{2}^{2}=\sum_{j=0}^{n-1} |x_{i}.\gamma_{j}|_{2}^{2}.
    \end{equation*}
    Let $x_{i}\in \mathbb{C}G_{i}$ and $\varepsilon >0 $ is given.  We apply the definition of norm operator  on operator $x_{i}$, then we  conclude that there is a $\gamma=\sum_{i=0}^{n-1}\gamma _{i}$ with $|\gamma|_{2}^{2}=\sum_{i=0}^{n}|\gamma_{i}|_{2}^{2}=1$ such that
   \begin{equation*}
   \sum_{j=0}^{n-1} |x_{i}.\gamma_{j}|_{2}^{2}=|x_{i}.\gamma|_{2}^{2}\geq (|x_{i}|_{op}-\varepsilon)^{2}
    \end{equation*}

     This shows that there is a $j\in Z_{n}$  such that
      \begin{equation*}
      \frac{(|x_{i}|_{op}-\varepsilon)^{2}}{n}\leq|x_{i}.\gamma _{j}|_{2}^{2}\leq|x.\gamma_{j}|^{2}\leq|x|_{op}^{2}   \end{equation*}

      since $|\gamma_{j}|_{2}\leq 1$.  So $|x|_{op}\geq \frac{|x_{i}|_{op}}{\sqrt{n}}$, for $i=0,1,\ldots n-1$.  We conclude  that the operator norm on $\mathbb{C}G$ is equivalent to direct sum norm
$|\sum x_{i}|_{op}=\sum |x_{i}|_{op} $.  Now if we put $E=C_{red}^{*}(G)$,  $F=\mathbb{C}G$,
$F_{i}=\mathbb{C}G_{i}$ and apply the above lemma to these spaces, the proof of the proposition would be completed
\end{proof}

\begin{remark1}Let $G$ and $H$ be two groups where $H$ is finite. Similar to above we can define an $H$-partition for $G$. It is a decomposition $G=\coprod_{h\in H} G_{h}$ such that $G_{h_{1}} G_{h_{2}} \subseteq G_{h_{1}h_{2}}$. This is equivalent to say that $G_{e}$ is a normal subgroup of $G$ whose quotient is isomorphic to $H$. In the same manner as above we can prove that an $H$-partition for a group $G$, gives an $H$- graded structure for $C_{red}^{*}(G)$. Moreover it is clear from the definition that an element $g\in G$ which is not in $G_{e}$, can be considered as a nontrivial homogenous element of $H$-graded algebra $C_{red}^{*}(G)$.

\end{remark1}
\; A part of the philosophy of noncommutative geometry is to translate  the classical facts about compact topological spaces into language of (noncommutative) Banach or $C^{*}$ algebras, see \cite {CONNES} or \cite {ENG}.  According to the Gelfand - Naimark theorem, there is a natural contravariant   functor from the category of compact Hausdorff topological space to the category of unital complex $C^{*}$ algebras. This functor assign to $X$, the commutative $C^{*}$ algebra $C(X)$ of all complex valued continuous functions on $X$.  It also assign to  a continuous function $f;X \rightarrow Y$, the C* algebra morphism $f^{*}:C(Y)\rightarrow C(X)$ by $f^{*}(\phi)=\phi \circ f $. Conversely every unital commutative $C^{*}$ algebra $A$ is isomorphic to $C(X)$ for some compact topological space $X$ and  every morphism from $C(Y)$ to $C(X)$ is equal to $f^{*}$ for  some continuous function $f:X \rightarrow Y$. In particular constant maps from $X$ to $Y$  correspond to morphism from $C(Y)$ to $C(X)$ with one dimensional range. Most of the statements about topological spaces have an algebraic analogy in the world of $C^{*}$ algebras. For example it can be easily shown that a topological space $X$ is connected if and only if the algebra $C(X)$ has no nontrivial idempotent, where a nontrivial idempotent  $a$ in an algebra $A$ is an element $a\in A$  such that $a^{2}=a$ and $a\neq 0, 1$.   So nonexistence of nontrivial idempotent for a (noncommutative) Banach algebra $A$ is interpreted as "noncommutative connectedness". \\
 For a topological space $X$, considering the above functor,  it is natural to identify $C(I\times X)$,  with $C(I,\;C(X))$, where $I$ is the unit interval. So in order to obtain a homotopy theory, for a Banach algebra $A$,  we define the Banach algebra  AI as follows:
 \begin{dfn}
  Let $A$ be a Banach algebra, we denote by $AI$, the Banach algebra  of all continuous  $\gamma :I \rightarrow A $ with the standard operations and norm. We define $\pi_{i} ;AI \rightarrow A$ with $\pi_{i}(\gamma)=\gamma(i)$, for $i=0,1$
  \end{dfn}
We say that two morphism $\alpha,\;\beta :A \rightarrow B$ are homotopic if there is a morphism $\Theta :A\rightarrow BI$ such that $\Theta_{0}=\alpha,\; \Theta_{1}=\beta$, where $\Theta_{i}=\pi_{i}\circ\Theta$. A morphism $\alpha  :A\rightarrow B$ is called  null homotopic if it is homotopic to a  morphism with one dimensional range. Obviously this is  a natural Banach algebraic analogy of classical null homotopicity. \\

For a unital Banach algebra $A$ with unit element 1, the element $\lambda.1$ is simply shown by $\lambda$  where $\lambda\in\mathbb{C}$ is  an  scalar.
For an element $a\in A$, we denote $\textit{sp}(a)$ for  all $\lambda\in \mathbb{C}$ such that $a-\lambda$ is not invertible. The spectral radius of $a$, which is denoted by $\textit{spr}(a)$, is defined as the infimum of all  $r$ such that the disc of radius $r$ around the origin contains $\textit{sp}(a)$.

\section{Main Result}
In the next theorem, which is the main result of this paper, $k$ is an arbitrary  positive integer, $G$ is a finite Abelian group and $A$ is a unital complex Banach algebra.
\begin{thm1}
Let A be a G-graded Banach algebra with no non-trivial idempotents. Let
$a\in A$ be a nontrivial homogenous element. Then 0 belongs to the convex
hull of the spectrum $sp( a^{k})$.\\
Further, if A is commutative, and a is invertible, then $a^{k}$ and 1 do not lie
in the same connected component of the space of invertible elements G(A).

\end{thm1}

\begin{proof}
 Without lose of generality we assume that  $G=Z_{n}$. So we have a multiplicative operator $T : A \rightarrow A$
 with $T^{n}=Id$. Let $a$ be a nontrivial homogenous element so  $T(a)=\lambda a$ where $\lambda\neq 1$ is a root of the unity. Assume  for contrary that the  convex hull of $sp(a^{k})$ does not contain $0$. Then $sp(a^{k})$ can be included in a disc with center $z_{0}$, for some $z_{0}\in\mathbb{C}$,  which is a subset of a branch of logarithm. Using holomorphic functional calculus as described in \cite[Chapter 10]{Rudin},  we can find a convergent series $\frac{}{}b=\sum_{i=0}^{\infty} c_{i}(a^{k}-z_{0})^{i}$ where $\exp b
=a^{k}$ where $c_{i}$'s are complex numbers.\\
 So $T(b)=\sum_{i=0}^{\infty} c_{i}(\lambda^{k}a^{k}-z_{0})^{i}$.
 Thus $T(b)$ and $b$ are power series in $a$. In particular $a b=b a $ and $T^{i}(b)T^{j}(b)=T^{j}(b)T^{i}(b)$, for all $i$ and $j$.  So we have $(a^{-1}\exp\frac{b}{k} )^{k}=1 $. Then $\textit{sp}( a^{-1}\exp\frac{b}{k})$ is a subset of the set of all $k$-th roots of unity. On the other hand since $A$ has no nontrivial idempotent , the spectrum of each element must be connected. So we  can assume that $\textit{sp}( a^{-1}\exp\frac{b}{k})=\{1\}$, otherwise we multiply $a$ with an appropriate  root of unity. Put $q= a^{-1}\exp\frac{b}{k}-1$. Then $q$ is  a quasinilpotent element of $A$, that is its spectral radius $\textit{spr}(q)=0$.  We have
\begin{equation}
 exp\frac{b}{k}=a+qa=a+q'
\end{equation}
where $q'$ is  a quasinilpotent too because $aq=qa$ and note that $\textit{spr}(aq)\leq \textit{spr}(a)\textit{spr}(q)$ for commuting elements $a$ and $q$, see \cite[page 302]{Rudin}.  Moreover we have

 \begin{equation}
 exp \frac{T(b)}{k}=T(a)+T(qa)=\lambda a+T(q')= \lambda a+q_{1}
\end{equation}

where $q_{1}$ is a quasinilpotent.\\

 On the other hand $(a+q)^{-1}-a^{-1}= -(a+q)^{-1}q a^{-1}$ then $(a+q)^{-1}=a^{-1}+q''$ for some quasinilpotent $q''$. Thus we have $\exp-\frac{b}{k}=a^{-1}+q''$. We obtain from (1) and (2) and the latest equation $\exp\frac{T(b)-b}{k}=\lambda+q_{2}$ for some quasinilpotent $q_{2}$. Then $\textit{sp}(T(b)-b)$ is a single point $\{\mu\}$ different from zero, since $\lambda\neq1$,  so $T(b)=b+\mu +y$ for some quasinilpotent $y$.
 Then $y=T(b)-b-\mu$ is a quasinilpotent element that can be expanded as a power series in $a$. Note that we emphasize on this power series expansion only for commutativity purpose.   We have $T(b)=b+\mu+y$ so by induction we obtain
 \begin{equation*}
  b=T^{n}(b)=b+n\mu+y+T(y)+T^{2}(y)+\ldots  +T^{n-1}(y).
   \end{equation*}
   Since $T$ is multiplicative and $y$ has a power series expansion in $a$, we conclude that $T^{j}(y)$'s are commuting quasinilpotent elements so their sum is quasinilpotent too, see \cite[page 302]{Rudin}.  This implies that $n\mu$ is  quasinilpotent which is a contradiction. This completes the first part of the theorem.\\
    Now assume that $A$ is commutative  and $a$ is a nontrivial homogenous element such that $a^{k}$ is in the same connected component as the identity. Then  there exist an  element $b\in A$ with $\exp b=a^{k}$. Obviously the same argument as above, but without needing to expansion of b as a power series in a, leads to a contradiction. So the  proof is complete.
\end{proof}

The following  corollaries are immediate consequence of the above theorem:
\begin{coro1}
Let $X$ be a compact locally path connected and simply connected space and $\phi : X\mapsto X$ be a  homeomorphism  of order $n$. Assume that  $\lambda\neq 1$ is an nth root of the unity. Then for every continuous function  $f:X\rightarrow  \mathbb{C}$, there is a point $x\in X$ such that $\sum_{i=0}^{p-1} \lambda^{i}f(\phi^{i}(x))=0$
\end{coro1}

 \begin{proof}\; Let $A$ be the commutative Banach algebra of   continuous functions $f:X\rightarrow  \mathbb{C}$ with the usual structures. Since $X$ is connected , $A$ has no nontrivial idempotent.  Define the continuous automorphism $T:A\mapsto A$ by $T(f)=f\circ\phi$. $T$ satisfies $T^{n}=1$, so we have a $Z_{n}$-graded structure for $A$ in the form $A=\bigoplus_{i=0}^{n-1} \ker (A-\lambda^{i})$.  Put $g=\sum_{i=0}^{p-1} \lambda^{i}f(\phi^{i}(x))$, then $T(g)=\lambda^{n-1}g$, so $g$ is a nontrivial homogenous element of $A$. If $g(x)\neq0$ for all $x\in X$ then $g$ is  an invertible element of $A$ which is not in the same connected component as the identity, by the above theorem.
On the other hand, consider the covering space $\exp  :\mathbb{C} \rightarrow \mathbb{C}-\{0\}$. Since $X$ is simply connected and locally path connected, there is a lifting $h$ of $f$, that is  $h\in A$ with $\exp h =g$, using the standard lifting lemma in the theory of covering space, see  \cite [proposition 1.33] {HATCHER}. So $g$ is  a logarithmic element and must be in the same connected component as the identity. This contradicts to the fact that $g$ can not be connected to the identity. This completes  the proof.
\end{proof}

Putting $X=S^{2}, \; \phi(x)=-x$ and $\lambda=-1$ we obtain the classical two dimensional Borsuk Ulam theorem  as follows:
\begin{coro2}
For a   continuous function $f: S^{2}\rightarrow \mathbb{C}$  there must exist a point $x\in S^{2}$ with $f(x)=f(-x)$

\end{coro2}

\begin{exa2}
Put $S^{2}=\{(x,y,z)\in \mathbb{R}^{3}:x^{2}+y^{2}+z^{2}=1\}$ and let $f:S^{2}\rightarrow\mathbb{R}^{2}$ be a continuous
function. Then there is a point $(x_{0},y_{0},z_{0}) \in S^{2}$ such that $f(y_{0},-x_{0},-z_{0})-f(-x_{0},-y_{0},z_{0})+f(-y_{0},x_{0},-z_{0})-f(x_{0},y_{0},z_{0})=0$. To prove this, consider the four order homeomorphism $\phi$ of $S^{2}$ with
$\phi(x,y,z)=(-y,x,-z)$. Now apply corollary 1 with $\lambda=-1$.
\end{exa2}

In the  following  proposition we  give  a  generalization of  corollary 2, in term of action of compact groups on certain topological  spaces:
\begin{pro3}
Let $G$  be  a  compact  Abelian group with Haar measure $\mu$ which  acts on a simply  connected  and locally path connected compact space $X$.
Assume that $\phi$ is  a  nontrivial character of $G$. Then for every  continuous  function $f\in C(X)$, there exists
$x_{0}\in X$  such that $\int_{G} \phi(g)f(g.x_{0})d\mu =0$

\end{pro3}

\begin{proof}
Define $F\in C(X)$  with $F(x)= \int_{G} \phi(g)f(g.x)d\mu $. Choose  $h\in G$ with $\lambda=\phi(h)\neq 1$. Define the automorphism  $T_{h}$ on $C(X)$  with $T_{h}(K)(x)=K(h^{-1}.x)$. Then $T_{h}(F)= \lambda F$. If $F$  is  a  non vanishing  function, similar to the  same  argument as in the proposition 2, there exists $k\in  C(X)$ such that $F=e^{k}$. Thus $e^{T_{h}(k)}=T_{h}(e^{k})=T_{h}(F)=\lambda F$. This  shows  that $e^{T_{h}(k)-k} =\lambda \neq 1$.  Since $X$ is  connected, there exists a  non zero $l\in \mathbb{C}$  such that $T_{h}(k)=k+l$ which implies that
$T_{h}^{n}(k)=k+nl$ for  all  $n\in  \mathbb{N}$. This is  a  contradiction because every automorphism of  a  $C^{*}$ algebra is  an isometry.

\end{proof}

 The following  corollary is  an obvious consequence of the last part of the main  theorem. In this corollary and its sequel, $G(A)_{0}$ is the connected component of the identity.\\

\begin{coro3}
Let $A$ be an idempotentless  commutative Banach algebra which is graded  by a finite Abelian group such that a nontrivial homogenous element is invertible. Then $\frac{G(A)}{G(A)_{0}}$ is an infinite group.
\end{coro3} \

\section{Further Questions and Remarks}
In this section we present some questions which naturally arise from the main theorem and the corollaries of the previous sections.\\
 First we discuss about a pure  algebraic analogy of the  corollary 3. For this purpose we need some elements of stable rank theory and $K$-theory  for both Banach algebras  and complex algebras.  We say that a commutative  Banach algebra $A$ has topological stable rank one if $G(A)$ is dense in $A$, see \cite {RIEFEL}. Let $A$ be a commutative complex algebra. We say that a pair $(a_{1}, a_{2})$ is invertible if  there is  a pair $(b_{1},b_{2})$ such that $a_{1}b_{1}+a_{2}b_{2}=1$ . We say $A$ has Bass   stable rank one if for every invertible pair $(a,b)\in A^{2}$, there is an element $x\in A$ such that $a+bx$ is invertible, see \cite {BASE}.  It is mentioned in   \cite{RIEFEL} that for a commutative C*  algebra, the Bass stable rank  coincide with the topological stable rank.\\
  For  a complex Banach algebra  $A$, let $GL_{n}(A)$ be the space of all n by n invertible matrices with entries in $A$. There is  a natural topology on this  space and there is a natural embedding of  $GL_{n}(A)$ into  $GL_{n+1}(A)$ which sends a matrices $B$ to diag(B,1), so we have an inductive system with  $GL_{n}(A)$'s.  Put $GL(A)=\varinjlim GL_{n}(A)$. Then $GL(A)$ has a natural inductive limit topology and algebraic structure . Define $K_{1}(A)=\frac{GL(A)}{GL(A)_{0}}$ where $GL(A)_{0}$ is the connected component of the identity. This Abelian group $K_{1}(A)$ is the standard $K_{1}$ functor defined on the category  of Banach algebras. For a complex algebra $A$,  there is  a pure algebraic $K_{1}^{alg}$ functor defined as the quotient of $GL(A)$ by its commutator. On the other hand it is well known that for an stable rank one Banach algebra $A$, $\frac{G(A)}{G(A)_{0}}$ is naturally isomorphic to $K_{1}(A)$, see \cite{HAAGERUP}\\
  So considering this isomorphism, the equality  of the Bass and topological stable rank for commutative $C^{*}$ algebras and the  above corollary,  it is natural to ask the next pure algebraic question:

\begin{question1}
Let   $A$  be an idempotentless involutive and commutative complex algebra with Bass stable rank one which is graded by a finite Abelian group. Assume that a nontrivial homogenous element of $A$ is invertible. Does this implies that $K_{1}^{alg}(A)$ is an infinite group?
\end{question1}

In the following example, we drop the commutativity of both grading
group G and the graded algebra A. We observe that in
this case the main theorem is no longer valid.

\begin{exa3}

Let $F_{2}$ be the free group on two generators $x$ and $y$.  We shall see that there is a normal subgroup $G_{0}$ of $F_{2}$ which does not contain $yxy^{-1}x^{-1}$ and its quotient is isomorphic to $S_{3}$, the permutation group on three  elements. Assuming the existence of such  subgroup $G_{0}$,  we obtain an $S_{3}$-graded structure for  $C_{red}^{*}(F_{2})$ for which $yxy^{-1}x^{-1}$, as an element of $C_{red}^{*}(F_{2})$,  is a nontrivial homogenous element.  On the other hand it is well known that  this algebra has no nontrivial idempotent and $yxy^{-1}x^{-1}\in C_{red}^{*}(F_{2})$,   is in the same  connected component as  the identity, see \cite{HAAGERUP}. To prove the existence of such subgroup $G_{0}$, we first note that $F_{4}$, the free group on four generators, can be considered as an index three subgroup of $F_{2}$ which does not contain $yxy^{-1}x^{-1}$. This can be proved using  certain covering space of figure-8 space as illustrated in \cite[p.58]{HATCHER}. Let $h:\widetilde{X} \rightarrow  X$ be  a covering space with $h(q)=p$. Then the induced map $h_{*}: \pi_{1} (\widetilde{X},q) \rightarrow \pi_{1} (X, p)$ is an injective map. $h_{*}(\pi_{1} (\widetilde{X},q))$ is isomorphic to $\pi_{1} (\widetilde{X},q)$ which index , as a subgroup of $\pi_{1}(X,p)$ is  equal to the cardinal of a fibre of the covering space.  Moreover a loop $\gamma \in \pi_{1}(X, p)$ lies in the rang of $h_{*}$ if and only if the unique  lifting $\widetilde{\gamma}$ of $\gamma$ with starting point $q$ is a loop with base point $q$,
see  \cite[proposition 1.31]{HATCHER}.\\
 Now consider the 3-fold  covering space which is illustrated in the figure below. The fundamental group of the total space is $F_{4}$ and the fundamental group of the base space is $F_{2}$. The loop $yxy^{-1}x^{-1}$ is not in the range of projecting map of the covering  because its lifting with base point $q$ ends to a different  point $O$
 \\[1cm]

 \includegraphics[width=9cm,height=10cm]{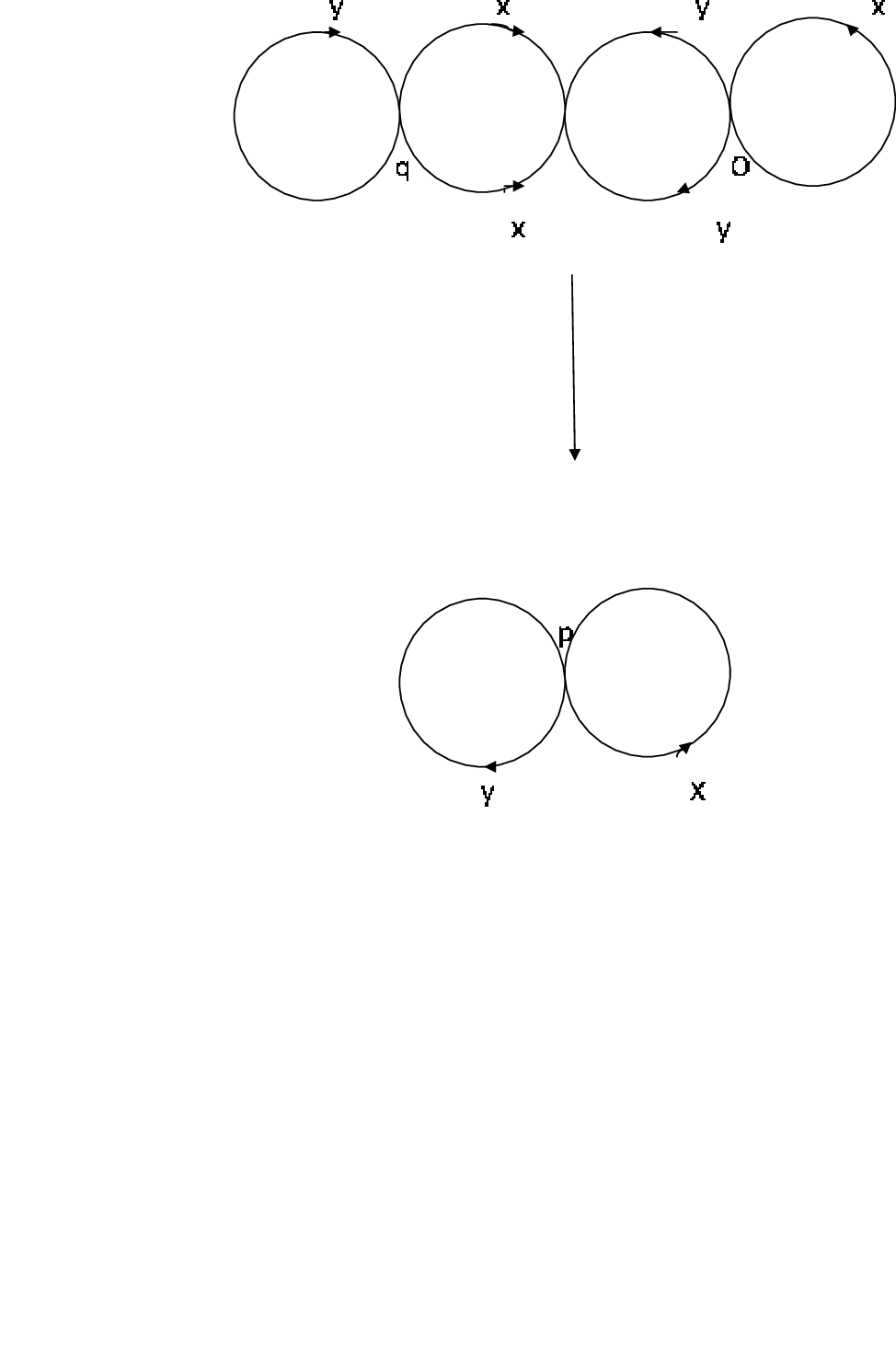}

  Put $G=F_{2}$ and $b=yxy^{-1}x^{-1}$. Then $G$ has a subgroup $H$ of index 3 which does not contain a commutator element $b$. We obtain a morphism $\beta: G  \rightarrow S_{3}$ with the standard action of G on the set of left cosets of $H$.
  Then $\ker \beta$  is  a normal subgroup of $G$ which is contained in $H$ and does not contain $b$. $G/\ker\beta$ is isomorph to a subgroup of $S_{3}$, on the other hand $G/\ker\beta$ is not Abelian since $\ker\beta$ does not contain a commutator element $b$. This shows that $\ker\beta$ is a normal subgroup of $G=F_{2}$ which quotient is isomorphic to $S_{3}$ and does not contain the commutator element $b=yxy^{-1}x^{-1}$.   So the above construction gives  an example of a non commutative Banach algebra $C_{red}^{*}(F_{2})$, without nontrivial idempotent,  which is graded by a non Abelian finite group $S_{3}$, such that a nontrivial homogenous element $yxy^{-1}x^{-1}$ lies in the same connected component as the identity.\\

\end{exa3}

But could we give any such example with a finite   Abelian group $G$?  In  other words, can we drop the hypothesis of commutativity of Banach algebra $A$ from the second part of the main theorem. This is a motivation for the next question which can be considered as a noncommutative analogy of the two dimensional Borsuk \verb|-|Ulam theorem:

\begin{question2}
Let $A$ be a  Banach algebra without nontrivial idempotent which is equipped with a $G$-graded structure, where $G$ is a finite Abelian group. Can one prove that the connected component of the identity has null intersection with nontrivial homogenous elements?  As a particular case,  put $A=C_{red}^{*}(F_{2})$, with the $Z_{2}$-graded structure correspond to the $Z_{2}$-partition of $F_{2}$  to the union of odd and even words. Can a linear combination of odd words be connected to the identity. Is the decomposition of this algebra to even and odd words, the only $Z_{2}$  graded structure for $A$ , up to graded isomorphism?
\end{question2}

There is an affirmative answer to a particular case of the second part of the above question. We thank professor
 A. Valette for his affirmative answer in this case. This affirmative answer is mentioned in the following proposition.

 \begin{pro4}
  Every linear combination of two odd words $a$ and $b$ in $C_{red}^{*}(F_{2})$ can not be connected to the identity
\end{pro4}

   \begin{proof}
  For a contrardiction assume that $ra+sb$  can be connected to the identity in the space of invertible elements, where $r$ and $s$ are two complex numbers. With an small perturbation we can assume that $|r|<|s|$.(Note that with this perturbation the connected component does not change). Now put $z_{0}=r/s$, so $z_{0}a+b$ lies in the same connected component as the identity with $|z_{0}|<1$. Then the curve $tz_{0}a+b,\; t\in [0 \;\;1] $ is a curve which lies in the space of invertible elements of $C_{red}^{*}(F_{2})$. Because $ba^{-1}$ is a unitary element  so its spectrum does not contain an element $tz_{0}$ where $t\in [0,\; 1]$. So $b$ as an odd element lies in the same connected component as the identity. But the only words which can be connected to the identity are members of the commutator subgroup of $F_{2}$. This contradiction shows that a linear combination of two odd words can not be connected to the identity
  \end{proof}

What is a Banach algebraic formulation of the higher dimensional Borsuk-Ulam theorem?
 In order to obtain a noncommutative version of this theorem  , we restate the classical case as follows;\\
  Let $f_{1},f_{2},\cdots, f_{n}$ be $n$ odd continuous real valued functions on $S^{n}$, then the function
 $ \sum_{i=1}^{n} f_{i}^{2}$ is not an invertible element of $C(S^{n})$ or equivalently is not in the same connected component as the identity (Sine every invertible element can be connected to the identity). In fact $f_{i}$'s are self adjoint elements of $C(S^{n})$ which are odd elements of the standard $Z_{2}$ graded structure of $C(S^{n})$.
 Now a relevant noncommutative version of this statement can be presented as the next question.  So it seems  natural to ask that for what type of noncommutative spheres the answer to the next question is affirmative?

\begin{question3}

 Assume that $A$ is a $Z_{2}$  graded   noncommutative $n$-sphere and  $a_{1}, a_{2},\ldots a_{n}$ are self adjoint   elements of $A$ which are odd elements of this graded algebra.  Is it true to say that $\sum_{i=1}^{n} a_{i}^{2} $ is either non invertible or is not in the same connected component as the identity?

\end{question3}

\begin{remark4}
A family of noncommutative spheres is a  family of C* algebras $A_{\theta}, \;\theta\in \mathbb{R}$ with some relations as a natural generalization of the algebra of continuous functions on n-sphere. For more information on noncommutative spheres, see \cite {NCS} or \cite{NAT}  \\
\end{remark4}

Another candidate for the noncommutative analogy of the Borsuk Ulam theorem can be presented as follows;

Consider the equivalent statement of the Borsuk Ulam theorem which says that an odd continuous map$f : S^{n} \rightarrow S^{n}$ is not null homotopic, namely it is not homotopic   to a constant map.   We try to translate this statement  into the language of Banach or $C^{*}$ algebras. The antipodal map $\phi(x)=-x$ defines  an order two automorphism $T ; C(S^{n})\rightarrow C(S^{n})$ with $T(g)=g \circ \phi$  which naturally gives a $Z_{2}$ graded structure for $C(S^{n})$. Similarly an odd map $f : S^{n} \rightarrow S^{n}$ defines a morphism $\alpha : C(S^{n})\rightarrow C(S^{n})$  which satisfies
$T \alpha =\alpha T$. This means that $\alpha$ is a graded morphism. So we ask: for what type of noncommutative spheres (spaces),  the answer to the next question is affirmative?

\begin{question4}
Let $A$ be a noncommutative sphere with a nontrivial $Z_{2}$ graded structure and $\alpha: A\rightarrow A$ be a graded morphism. Is it  true to say that $\alpha$ is not a null homotopic morphism ?
\end{question4}.

\begin{remark5}

We explain that the main theorem gives us a weaker version of the Kaplansky-Kadison conjecture. This conjecture says that for a torsion free group $\Gamma$,  $C_{red}^{*}(\Gamma)$ has no nontrivial idempotent, see \cite{VALETTE}.  Now as a consequence of our theorem  we present the following question as a weaker version of the Kaplansky-Kadison conjecture:

\begin{question5}
Let $\Gamma$ be a torsion free group and $C_{red}^{*}(\Gamma)$ is equipped  with a $Z_{n}$- graded structure. Is it possible that the convex hull of the spectrum of a nontrivial homogenous element does not contain the origin?
\end{question5}

\end{remark5}

\textbf{Acknowledgment.}\;
The author would like to thank the referee for very valuable suggestions.\\
\bigskip\\

\end{document}